\documentclass[12pt,draft]{article}

\usepackage[cp1251]{inputenc}
\usepackage[T2A]{fontenc}
\usepackage[ukrainian,english,russian]{babel}
\usepackage{amsmath,amsfonts,amssymb}
\usepackage{geometry}
\begin{document}

\begin{flushleft}


\vspace{0,5cm}

\large\textbf{Татьяна Н. Зинченко, Александр А. Мурач }\\
\normalsize(Институт математики НАН Украины, Киев, Украина)

\vspace{0,5cm}

\large\textbf{ЭЛЛИПТИЧЕСКИЕ ПО ДУГЛИСУ--НИРЕНБЕРГУ СИСТЕМЫ \\ В ПРОСТРАНСТВАХ
ХЕРМАНДЕРА}

\vspace{0,5cm}

\textbf{Tatjana N. Zinchenko, Aleksandr A. Murach}\\
\normalsize(Institute of Mathematics of NAS of Ukraine,  Kyiv, Ukraine)

\vspace{0,5cm}

\large\textbf{DOUGLIS--NIRENBERG ELLIPTIC SYSTEMS \\
IN H\"ORMANDER SPACES}

\end{flushleft}

\vspace{0,5cm}

\noindent Исследуются равномерно эллиптические в $\mathbb{R}^{n}$ по
Дуглису--Ниренбергу системы в классе гильбертовых пространств Хермандера. Последние
параметризуются радиальным функциональным параметром, который RO-меняется на
$+\infty$ как функция от $(1+|\xi|^{2})^{1/2}$, где $\xi\in\mathbb{R}^{n}$. Доказана
априорная оценка решений и исследована их регулярность. Получено достаточное условие
нетеровости этих систем.

\vspace{0,5cm}

\noindent We investigate Douglis--Nirenberg uniformly elliptic systems in
$\mathbb{R}^{n}$ on a class of H\"ormander inner product spaces. They are
parametrized with a radial function parameter which is RO-varying at $+\infty$,
considered as a function of $(1+|\xi|^{2})^{1/2}$ with $\xi\in\mathbb{R}^{n}$. An
a'priori estimate for solutions is proved, and their interior regularity is studied.
A sufficient condition for the systems to have the Fredholm property is given.

\vspace{0,5cm}

\noindent\textbf{1. Введение.} Общие эллиптические системы дифференциальных
уравнений смешанного порядка были введены А.~Дуглисом и Л.~Ниренбергом
\cite{DouglisNirenberg55}. Содержательные примеры таких систем встречаются в
гидродинамике и теории упругости. Эллиптические по Дуглису--Ниренбергу системы
возникают также при сведении скалярных эллиптических уравнений к системам уравнений
первого порядка и при сведении эллиптических краевых задач на край многообразия
\cite{Agranovich94, Agranovich97, WlokaRowleyLawruk95}.

Эллиптические уравнения и системы обладают рядом характерных свойств в шкалах
пространств Гельдера--Зигмунда и Соболева: априорные оценки решений, повышение
гладкости решений, нетеровость эллиптических операторов. Указанные свойства имеют
важные приложения в теории эллиптических краевых задач, в теории индекса
эллиптических операторов, в спектральной теории дифференциальных операторов и ряд
других (см. обзоры \cite{Agranovich94, Agranovich97} и приведенную там литературу).

В этой связи представляет интерес исследование эллиптических уравнений и систем в
различных классах функциональных пространств, характеризующих свойства регулярности
функций/распределений более тонко, чем классические шкалы Гельдера--Зигмунда и
Соболева. Для этого естественно использовать пространства, в которых показателем
гладкости служит не числовой, а функциональный параметр. Широкие классы таких
пространств были введены и систематически исследованы Л. Хермандером
\cite{Hermander63} (п. 2.2), и применены к изучению локальной регулярности решений
линейных дифференциальных уравнений, а также их систем (последние~--- в случае
постоянных коэффициентов) \cite{Hermander63, Hermander83}. В настоящее время
пространства Хермандера и их различные аналоги, именуемые пространствами обобщенной
гладкости, активно исследуются как сами по себе, так и с точки зрения приложений
[7--12].

Для приложений, особенно в спектральной теории, наиболее важным является случай
гильбертовых пространств. Среди них особый интерес представляют пространства,
интерполяционные относительно гильбертовой соболевской шкалы. Поскольку при
интерполяции наследуется ограниченность линейных операторов, а также их нетеровость
(при неизменном дефекте), то классы этих пространств служат удобным инструментом для
изучения свойств эллиптических уравнений и систем.

В настоящей статье исследуются равномерно эллиптические в $\mathbb{R}^{n}$ по
Дуглису--Ни\-рен\-бер\-гу системы в классе гильбертовых пространств Хермандера
$$
H^{\varphi}:=B_{2,\varphi(\langle\cdot\rangle)}=\{w\;\mbox{---
распределение}:\,\varphi(\langle\xi\rangle)(\mathcal{F}w)(\xi)\in
L_{2}(\mathbb{R}^{n},d\xi)\}, \eqno(*)
$$
где $\varphi$~--- произвольная RO-меняющаяся функция скалярного аргумента. Выбор
этого класса обусловлен тем, что он совпадает (с точностью до эквивалентности норм)
с классом всех гильбертовых пространств, интерполяционных относительно гильбертовой
соболевской шкалы \cite{MikhailetsMurach10} (п. 2.4.2),
\cite{MikhailetsMurach1106arxiv}. В статье установлены теоремы об априорных оценках
и о регулярности решений эллиптических систем в пространствах $(*)$, а также (при
дополнительных предположениях) теорема о нетеровости матричного эллиптического
оператора. Отдельный случай систем, равномерно эллиптических по Петровскому,
рассмотрен ранее в \cite{09UMJ3}.

Отметим, что для более узкого класса пространств Хермандера (уточненная соболевская
шкала) эллиптические уравнения и эллиптические краевые задачи исследованы
В.~А.~Михайлецом и вторым автором в серии работ, среди которых упомянем статьи
[15--20] и монографию \cite{MikhailetsMurach10}.

\textbf{2. Постановка задачи.} Пусть $p,n\in\mathbb{N}$. В евклидовом пространстве
$\mathbb{R}^n$ рассматривается система линейных дифференциальных уравнений
\begin{equation}\label{f1}
\sum _{k=1}^p A_{j,k}(x,D)u_{k}(x)=f_{j}(x),\quad j=1,\ldots,p,
\end{equation}
где
$$
A_{j,k}(x,D):=\sum_{|\mu|\leq r_{j,k}}a_{\mu}^{j,k}(x)D^\mu,\quad j,k=1,\ldots,p.
$$
Здесь и далее $\mu=(\mu_1,\ldots,\mu_n)$~--- мультииндекс с неотрицательными целыми
компонентами, $|\mu|$:= $\mu_1+\ldots+\mu_n$, $D_{j}:=i\partial/\partial x_{j}$,
$D^\mu(=D^\mu_x):=D_{1}^{\mu_1}\ldots D_{n}^{\mu_n}$, где $i$~--- мнимая единица, а
$x=(x_{1},\ldots,x_{n})\in\mathbb{R}^n$. Преобразование Фурье $\mathcal{F}$
переводит дифференциальный оператор $D^\mu$ в оператор умножения на функцию
$\xi^\mu:=\xi_{1}^{\mu_1}\ldots\xi_{n}^{\mu_n}$ аргумента
$\xi=(\xi_{1},\ldots,\xi_{n})\in\mathbb{R}^n$, двойственного к $x$.

Предполагается, что в системе $\eqref{f1}$ все коэффициенты $a_{\mu}^{j,k}(x)$~---
комплекснозначные функции, бесконечно дифференцируемые и ограниченные вместе со
всеми частными производными в $\mathbb{R}^n$. Класс таких функций обозначаем через
$C^{\infty}_{\mathrm{b}}$.

Решение системы $\eqref{f1}$ понимается в смысле теории распределений. Запишем ее в
матричной форме $Au=f$. Здесь $A:=(A_{j,k}(x,D))_{j,k=1}^{p}$~--- матричный
дифференциальный оператор, а $u=\mathrm{col}\,(u_{1},\ldots,u_{p})$,
$f=\mathrm{col}\,(f_{1},\ldots,f_{p})$~--- функциональные столбцы.

Предполагается, что система \eqref{f1} равномерно эллиптическая в $\mathbb{R}^{n}$
по Дуглису--Ни\-рен\-бер\-гу \cite{Agranovich94} (п. 3.2.b), т.~е. существуют наборы
вещественных чисел ${l_1,\ldots,l_p}$ и ${m_1,\ldots,m_p}$ такие, что:

i) $r_{j,k}\leq l_j+m_k$ для всех $j,k = 1,\ldots,p$;

ii) найдется число $c>0$, при котором
$$
|\det(A_{j, k}^{(0)}(x,\xi))_{j,k=1}^p|\geq c \quad\mbox{для любых}\quad x,
\xi\in\mathbb{R}^{n},\;|\xi| = 1.
$$
Здесь
$$
A_{j, k}^{(0)}(x,\xi):=\sum_{|\mu|= l_j+m_k}a_{\mu}^{j,k}(x)\,\xi^\mu
$$
--- главный символ дифференциального оператора $A_{j,k}(x,D)$ в случае
$r_{j,k}=l_j+m_k$, либо $A_{j, k}^{(0)}(x,\xi)\equiv0$ в случае $r_{j,k}<l_j+m_k$.

В отдельном случае, когда все числа $l_{j}=0$, система \eqref{f1} называется
равномерно эллиптической по Петровскому. Если, кроме того, все числа $m_k$ равны, то
она является равномерно эллиптической в обычном смысле.

Введем функциональные пространства, в которых исследуется система \eqref{f1}.

Пусть RO~--- множество всех измеримых по Борелю функций
$\varphi:[1,\infty)\rightarrow(0,\infty)$, для которых существуют числа $a>1$ и
$c\geq1$ такие, что
$$
c^{-1}\leq\frac{\varphi(\lambda t)}{\varphi(t)}\leq c\quad\mbox{для любых}\quad
t\geq1,\;\lambda\in[1,a]
$$
(постоянные а и с зависят от $\varphi\in\mathrm{RO}$). Такие функции называют RO
(или OR)-ме\-ня\-ю\-щимися на бесконечности. Класс RO-меняющихся функций введен
В.~Г.~Авакумовичем в 1936 г. и достаточно полно изучен (см. \cite{Seneta76}
(приложение~1), \cite{BinghamGoldieTeugels89} (п. 2.0~-- 2.2)).

Пусть $\varphi\in\mathrm{RO}$. Обозначим через $H^{\varphi}$ линейное пространство
всех распределении $w\in\mathcal{S}'$ таких, что их преобразование Фурье
$\widehat{w}:=\mathcal{F}w$ локально суммируемо по Лебегу в $\mathbb{R}^{n}$ и
удовлетворяет условию
$$
\int\limits_{\mathbb{R}^{n}}\varphi^2(\langle\xi\rangle)\,|\widehat{w}(\xi)|^2\,d\xi
<\infty.
$$
Здесь, как обычно, $\mathcal{S}'$~--- линейное топологическое пространство Шварца
медленно растущих комплекснозначных распределений, заданных в $\mathbb{R}^{n}$, а
$\langle\xi\rangle:=(1+|\xi|^{2})^{1/2}$~--- сглаженный модуль вектора
$\xi\in\mathbb{R}^{n}$. С точки зрения приложений к дифференциальным уравнениям нам
удобно трактовать распределения как \emph{анти}линейные функционалы на пространстве
$\mathcal{S}$ основных функций.

В пространстве $H^{\varphi}$ определено скалярное произведение распределений $w_1$,
$w_2$ по формуле
$$
(w_1,w_2)_{\varphi}:=\int\limits_{\mathbb{R}^{n}}\varphi^2(\langle\xi\rangle)\,
\widehat{w_1}(\xi)\,\overline{\widehat{w_2}(\xi)}\,d\xi.
$$
Оно задает на $H^{\varphi}$ структуру гильбертова пространства и определяет норму
$\|w\|_\varphi:=(w,w)_\varphi^{1/2}$. Это пространство сепарабельно; в нем плотно
множество $C_{0}^{\infty}$ бесконечно дифференцируемых функций на $\mathbb{R}^{n}$,
у которых носитель компактен.

Пространство $H^{\varphi}$~--- гильбертов изотропный случай пространств $B_{p,k}$,
введенных и систематически исследованных Л.~Хермандером \cite{Hermander63} (п. 2.2)
(см. также \cite{Hermander83} (п. 10.1)). Именно, $H^{\varphi}=B_{p,k}$, если $p=2$
и $k(\xi)=\varphi(\langle\xi\rangle)$ при $\xi\in\mathbb{R}^{n}$. Отметим, что при
$p=2$ пространства Хермандера совпадают с пространствами, введенными и изученными
Л.~Р.~Волевичем и Б.~П.~Панеяхом \cite{VolevichPaneah65} (\S~2).

Если $\varphi(t)=t^{s}$ для всех $t\geq1$ при некотором $s\in\mathbb{R}$, то
$H^{\varphi}=:H^{(s)}$ есть (гильбертово) пространство Соболева порядка $s$.

Отметим, что пространства $H^{\varphi}$ и $H^{1/\varphi}$ взаимно двойственны
относительно расширения по непрерывности полуторалинейной формы
$$
(w_{1},w_{2})_{\mathbb{R}^{n}}:=
\int\limits_{\mathbb{R}^{n}}\,w_{1}(x)\,\overline{w_{2}(x)}\,dx.
$$
(очевидно, $\varphi\in\mathrm{RO}\Leftrightarrow1/\varphi\in\mathrm{RO}$). Это
расширение обозначаем также через $(\cdot,\cdot)_{\mathbb{R}^{n}}$, а для
вектор-функций $u$ и $f$ полагаем
$(u,f)_{\mathbb{R}^{n}}:=(u_{1},f_{1})_{\mathbb{R}^{n}}+
\ldots+(u_{p},f_{p})_{\mathbb{R}^{n}}$, если слагаемые определены.

Обозначим $\varrho(t):=t$ при $t\geq1$. Матричный дифференциальный оператор $A$
является ограниченным оператором (см. ниже п.~4)
\begin{equation}\label{f2}
A:\,\bigoplus_{k=1}^p\,H^{\varphi\rho^{m_k}}\rightarrow
\bigoplus_{j=1}^p\,H^{\varphi\rho^{-l_j}}\quad\mbox{для каждого}
\quad\varphi\in\mathrm{RO}.
\end{equation}
Здесь $\varphi\rho^{m_k},\varphi\rho^{-l_j}\in\mathrm{RO}$ и поэтому определены
пространства Хермандера, фигурирующие в~\eqref{f2}.

В работе исследуются свойства оператора~\eqref{f2}.

\textbf{3. Основные результаты}. Сформулируем основные результаты статьи; их
доказательство будет дано ниже в п.~5.

\textbf{Теорема 1.} \it Пусть заданы функция $\varphi\in\mathrm{RO}$ и число
$\sigma>0$. Тогда существует число $c=c(\varphi,\sigma)>0$ такое, что для
произвольных вектор-функций
\begin{equation}\label{f3}
u\in\bigoplus_{k=1}^p\,H^{\varphi\rho^{m_k}},\quad
f\in\bigoplus_{j=1}^p\,H^{\varphi\rho^{-l_j}},
\end{equation}
удовлетворяющих уравнению $Au=f$ в $\mathbb{R}^{n}$, справедлива априорная оценка
\begin{equation}\label{f4}
\biggl(\,\sum_{k=1}^{p}\|u_{k}\|_{{\varphi\rho}^{m_k}}^2\biggr)^{1/2}\leq
c\,\biggl(\,\sum_{j=1}^{p}\|f_{j}\|_{{\varphi\rho}^{-l_j}}^2\biggr)^{1/2}+
c\,\biggl(\,\sum_{k=1}^{p}\|u_{k}\|_{{\varphi \rho}^{m_k-\sigma}}^2\biggr)^{1/2}.
\end{equation} \rm

Пусть $V$~--- произвольное открытое непустое подмножество пространства
$\mathbb{R}^{n}$. Исследуем внутреннюю регулярности решения эллиптической системы
$Au=f$ на $V$ в классе пространств Хермандера.

Обозначим
\begin{gather}\notag
H^{-\infty}:=\bigcup_{s\in\mathbb{R}}H^{(s)}=
\bigcup_{\varphi\in\mathrm{RO}}H^{\varphi},\quad
H^{\infty}:=\bigcap_{s\in\mathbb{R}}H^{(s)}=
\bigcap_{\varphi\in\mathrm{RO}}H^{\varphi}
\end{gather}
Эти определения корректны, как будет показано в п.~4. В пространствах $H^{-\infty}$
и $H^{\infty}$ вводятся топологии соответственно индуктивного и проективного
пределов. Отметим, что $H^{\infty}\subset C_{\mathrm{b}}^{\infty}$ в силу теоремы
вложения Соболева. Положим
\begin{gather}\notag
H^{\varphi}_{\mathrm{int}}(V):=\bigl\{w\in H^{-\infty}:\,\chi\,w\in H^{\varphi}\\
\mbox{для всех}\;\;\chi\in C^{\infty}_{\mathrm{b}}\;\;\mbox{таких,
что}\;\;\mathrm{supp}\,\chi\subset V,\;\;\mathrm{dist}(\mathrm{supp}\,\chi,\partial
V)>0\bigr\};\label{f5}
\end{gather}
здесь $\varphi\in\mathrm{RO}$. Топология в пространстве
$H^{\varphi}_{\mathrm{int}}(V)$ задается полунормами $w\rightarrow\|\chi
w\|_\varphi$, где функции $\chi$ те же, что и в $\eqref {f5}$. Если
$V=\mathbb{R}^{n}$, то $H^{\varphi}_{\mathrm{int}}(V)=H^{\varphi}$.

\textbf{Теорема 2.} \it Пусть $\varphi\in\mathrm{RO}$. Предположим, что
вектор-функция $u\in(H^{-\infty})^p$ является решением уравнения $Au=f$ на открытом
множестве $V\subseteq\mathbb{R}^{n}$, где $f_{j}\in
H_{\mathrm{int}}^{\varphi\rho^{-l_j}}(V)$ для всех $j=1,\ldots,p$. Тогда $u_{k}\in
H_{\mathrm{int}}^{\varphi\rho^{m_k}}(V)$ для всех $k=1,\ldots,p$. \rm

Отметим, что следует различать \emph{внутреннюю} и \emph{локальную} регулярность на
открытом множестве $V\subset\mathbb{R}^{n}$. Пространство распределений, имеющих
характеризуемую параметром $\varphi\in\mathrm{RO}$ локальную регулярность на этом
множестве, определяется следующим образом:
$$
H^{\varphi}_{\mathrm{loc}}(V):=\bigl\{w\in H^{-\infty}:\,\chi\,w\in
H^{\varphi}\;\;\mbox{для всех}\;\;\chi\in C^{\infty}_{0}\;\;\mbox{таких,
что}\;\;\mathrm{supp}\,\chi\subset V\bigr\}.
$$
В случае, когда множество $V$ ограничено, пространства
$H^{\varphi}_{\mathrm{int}}(V)$ и $H^{\varphi}_{\mathrm{loc}}(V)$ совпадают. Если же
$V$ не ограничено, то может быть строгое включение
$H^{\varphi}_{\mathrm{int}}(V)\subset H^{\varphi}_{\mathrm{loc}}(V)$. Для локальной
гладкости справедлив аналог теоремы~2; в ее формулировке следует лишь заменить
$\mathrm{int}$ на $\mathrm{loc}$ в обозначениях пространств. Он тривиально вытекает
из теоремы~2.

В качестве приложения этой теоремы имеем следующее достаточное условие непрерывности
частных производных решения $u$.

\textbf{Теорема 3.} \it Пусть заданы целые числа $k\in\{1,\ldots,p\}$,
$\lambda\geq0$ и функция $\varphi\in\mathrm{RO}$, удовлетворяющая условию
\begin{equation}\label{f6}
\int\limits_{1}^{\infty}\,t^{2\lambda+n-1-2m_k}\,\varphi^{-2}(t)\,dt<\infty.
\end{equation}
Предположим, что вектор-функция $u\in(H^{-\infty})^p$ является решением уравнения
$Au=\nobreak f$ на открытом множестве $V\subseteq \mathbb{R}^{n}$, где $f_j\in
H_{\mathrm{int}}^{\varphi\rho^{-l_j}}(V)$ для всех $j=1,\ldots,p$. Тогда компонента
$u_k$ решения имеет на множестве $V$ непрерывные частные производные до порядка
$\lambda$ включительно, причем эти производные ограничены на каждом множестве
$V_{0}\subset V$ таком, что $\mathrm{dist}(V_{0},\partial V)>0$. В частности, если
$V=\mathbb{R}^{n}$, то $u_k\in C_{\mathrm{b}}^{\lambda}$. \rm

Здесь $C_{\mathrm{b}}^{\lambda}$~--- банахово пространство всех функций
$w:\mathbb{R}^n\rightarrow\mathbb{C}$, имеющих непрерывные и ограниченные
производные в $\mathbb{R}^n$ порядка $\leq\lambda$.

Отметим, что аналоги теорем 1 -- 3 справедливы и для системы $A^{+}u=f$, формально
сопряженной к системе \eqref{f1}, поскольку обе они равномерно эллиптичны в
$\mathbb{R}^{n}$ (по Дуглису--Ниренбергу). Здесь, напомним,
$A^{+}:=(A_{k,j}^{+}(x,D))_{j,k=1}^{p}$, где
$$
A_{k,j}^{+}(x,D)u_{k}(x):=\sum_{|\mu|\leq
r_{k,j}}D^{\mu}({\overline{a_{\mu}^{k,j}(x)}}u_{k}(x));
$$
так что $(A^{+}u,v)_{\mathbb{R}^n}=(u,Av)_{\mathbb{R}^n}$ для произвольных
вектор-функций $u,v\in(\mathcal{S})^p$.

Системе $A^{+}u=f$ соответствует ограниченный оператор
\begin{equation}\label{f6a}
A^{+}:\,\bigoplus_{k=1}^p\,H^{\varphi\rho^{l_k}}\rightarrow
\bigoplus_{j=1}^p\,H^{\varphi\rho^{-m_j}}\quad\mbox{для каждого}
\quad\varphi\in\mathrm{RO}.
\end{equation}
Он сопряжен к оператору \eqref{f2}, где пишем $1/\varphi$ вместо $\varphi$,
относительно полуторалинейной формы $(\cdot,\cdot)_{\mathbb{R}^n}$.

Согласно теореме 2 ядра операторов \eqref{f2} и \eqref{f6a} совпадают с
пространствами
$$
N:=\{u\in(H^{\infty})^p:\,Au=0\},\quad N^{+}:=\{v\in(H^{\infty})^p:\,A^{+}v=0\}
$$
соответственно и не зависят от $\varphi$.

Следующие условия являются достаточными для нетеровости оператора \eqref{f2} (а
также оператора \eqref{f6a}):

а) $D^\alpha a_{\mu}^{j,k}(x)\rightarrow0$ при $|x|\to\infty$ для каждого
мультииндекса $\alpha$ с $|\alpha|\geq 1$, произвольных индексов
$j,k\in\{1,\ldots,p\}$ и мультииндекса $\mu$ с $|\mu|\leq r_{j,k}$;

б) существуют числа $c_{1}>0$ и $c_{2}\geq0$ такие, что
$$
|\det(A_{j, k}(x,\xi))_{j,k=1}^p|\geq c_{1}\langle\xi\rangle^{q}\quad\mbox{для
любых}\quad x,\xi\in\mathbb{R}^{n},\;|x|+|\xi|\geq c_{2}.
$$
Здесь $q:=l_1+\ldots+l_p+m_1+\ldots+m_p$, а
$$
A_{j, k}(x,\xi):=\sum_{|\mu|\leq l_j+m_k}a_{\mu}^{j,k}(x)\,\xi^\mu
$$
---~полный символ дифференциального оператора $A_{j,k}(x,D)$.

Напомним, что линейный ограниченный оператор $T:E_1\rightarrow E_2$, где $E_1$ и
$E_2$~--- банаховы пространства, называется нетеровым, если его ядро $\ker T$ и
коядро $\mathrm{coker}\,T:=E_2/T(X)$ конечномерны. У нетерового оператора $T$
область значений $T(X)$ замкнута в $E_2$, а индекс $\mathrm{ind}\,T:=\dim\ker
T-\dim\mathrm{coker}\,T$ конечен.

\textbf{Теорема 4.} \it Пусть выполняются условия \rm а) \it и \rm б)\it. Тогда для
каждого $\varphi\in\mathrm{RO}$ оператор $\eqref{f2}$ нетеров. Его область значений
совпадает с пространством
\begin{equation}\label{f7}
\biggl\{f\in\bigoplus_{j=1}^{p}H^{\varphi\rho^{-l_j}}: \,(f,v)_{\mathbb{R}^{n}}=0
\;\;\mbox{для всех}\;\;v\in N^{+}\biggr\},
\end{equation}
а индекс равен $\dim N-\dim N^{+}$ и не зависит от $\varphi$. \rm

Отметим, что условие б) влечет за собой условие ii) из определения равномерной
эллиптичности по Дуглису--Ниренбергу. В свою очередь, если предположить а), то
условие б) следует из нетеровости оператора $\eqref{f2}$ в соболевском случае
$\varphi=\varrho^{s}$ при хотя бы одном значении $s\in\mathbb{R}$ \cite{Rabier12}
(теорема 4.2).

Предположим, что выполняются условия а) и б). В случае, когда пространства $N$ и
$N^{+}$ тривиальны, оператор \eqref{f2} является гомеоморфизмом в силу теоремы~4 и
теоремы Банаха об обратном операторе. В общей ситуации гомеоморфизм удобно задавать
с помощью следующих проекторов.

Пусть $\varphi\in\mathrm{RO}$. Разложим пространства, в которых действует нетеров
оператор \eqref{f2}, в прямые суммы (замкнутых) подпространств:
\begin{gather*}
\bigoplus_{k=1}^{p}H^{\varphi\rho^{m_k}}=N\dotplus
\biggl\{u\in\bigoplus_{k=1}^{p}H^{\varphi\rho^{m_k}}:\,(u,w)_{\mathbb{R}^{n}}=0
\;\;\mbox{для всех}\;\;w\in N\biggr\},\\
\bigoplus_{j=1}^{p}H^{\varphi\rho^{-l_j}}=N^{+}\dotplus\eqref{f7}.
\end{gather*}
Такие разложения существуют, поскольку в них слагаемые имеют тривиальное
пересечение, и конечная размерность первого из них равна коразмерности второго.
Последнее следует из того, что в первой сумме факторпространство пространства
$\bigoplus_{k=1}^{p}H^{\varphi\rho^{m_k}}$ по второму слагаемому является
двойственным пространством к подпространству $N$ пространства
$\bigoplus_{k=1}^{p}H^{1/(\varphi\rho^{m_k})}$ (двойственность понимается
относительно формы $(\cdot,\cdot)_{\mathbb{R}^{n}}$). Аналогично и для второй суммы.

Обозначим через $P$ и $P^{+}$ соответственно (косые) проекторы пространств
$$
\bigoplus_{k=1}^p\,H^{\varphi\rho^{m_k}}\quad\mbox{и}\quad
\bigoplus_{j=1}^p\,H^{\varphi\rho^{-l_j}}
$$
на вторые слагаемые в указанных суммах параллельно первым слагаемым. Эти проекторы
(как отображения) не зависят от $\varphi$.

Тогда в силу теоремы~4 сужение оператора \eqref{f2} на подпространство
$P(\bigoplus_{k=1}^{p}H^{\varphi\rho^{m_k}})$ является гомеоморфизмом
$$
A:\,P\Bigl(\bigoplus_{k=1}^{p}H^{\varphi\rho^{m_k}}\Bigr)\leftrightarrow
P^{+}\Bigl(\bigoplus_{j=1}^p\,H^{\varphi\rho^{-l_j}}\Bigr).
$$

Аналогичный результат верен и для оператора \eqref{f6a}. Отметим, что его
нетеровость следует из теоремы~4, поскольку он сопряжен к нетеровому оператору
\eqref{f2} с параметром $1/\varphi$ вместо $\varphi$.

\textbf{4. Вспомогательные результаты.} Приведем некоторые полезные нам факты; они
будут использованы в доказательствах теорем 1~-- 4.

Отметим следующие свойства функционального класса RO (см., например,
\cite{Seneta76}, приложение~1, теоремы 1 и 2):

i) $\varphi\in\mathrm{RO}$ тогда и только тогда, когда
$$
\varphi(t)=\exp\biggl(\beta(t)+
\int\limits_{1}^{\:t}\frac{\alpha(\tau)}{\tau}\;d\tau\biggr)\quad \mbox{при}\quad
t\geq1,
$$
где вещественные функции $\alpha$ и $\beta$  измеримы по Борелю и ограничены на
полуоси $[1,\infty)$;

ii) для любой функции $\varphi\in\mathrm{RO}$ существуют числа
$s_{0},s_{1}\in\mathbb{R}$, $s_{0}\leq s_{1}$, и $c_{1}\geq1$ такие, что
\begin{equation}\label{f8}
c_{1}^{-1}\lambda^{s_{0}}\leq\frac{\varphi(\lambda t)}{\varphi (t)}\leq
c_{1}\lambda^{s_{1}} \quad\mbox{при}\quad t\geq1,\;\;\lambda\geq1.
\end{equation}

Для функции $\varphi\in\mathrm{RO}$ определены и конечны нижний и верхний индексы
Матушевской \cite{BinghamGoldieTeugels89} (п. 2.1.2):
\begin{gather*}
\sigma_{0}(\varphi):=
\sup\,\{s_{0}\in\mathbb{R}:\,\mbox{верно левое неравенство в \eqref{f8}}\},\\
\sigma_{1}(\varphi):=\inf\,\{s_{1}\in\mathbb{R}:\,\mbox{верно правое неравенство в
\eqref{f8}}\}.
\end{gather*}
Из формулы $\eqref{f8}$ при $t=1$ следуют непрерывные и плотные вложения
\begin{equation}\label{f9}
H^{(s_1)}\hookrightarrow H^{\varphi}\hookrightarrow H^{(s_0)}\quad\mbox{для всех
чисел}\quad s_{1}>\sigma_{1}(\varphi),\;\;s_{0}<\sigma_{0}(\varphi).
\end{equation}
Отсюда вытекает корректность определения пространств $H^{-\infty}$ и $H^{\infty}$,
данного в п.~3.

Пространство $H^{\varphi}$, фигурирующее в \eqref{f9}, есть результат интерполяции с
подходящим функциональным параметром пары соболевских пространств $H^{(s_0)}$ и
$H^{(s_1)}$. Напомним определение этой интерполяции в случае общих гильбертовых
пространств и некоторые ее свойства \cite{MikhailetsMurach10} (п. 1.1, 2.4.2). Для
наших целей достаточно ограничиться сепарабельными пространствами.

Пусть задана упорядоченная пара $X:=[X_{0},X_{1}]$ сепарабельных комплексных
гильбертовых пространств $X_{0}$ и $X_{1}$ такая, что выполняется непрерывное и
плотное вложение $X_{1}\hookrightarrow X_{0}$. Пару $X$ называем допустимой. Для нее
существует изометрический изоморфизм $J:X_{1}\leftrightarrow X_{\,0}$ такой, что
$J$~--- самосопряженный положительно определенный оператор в пространстве $X_{0}$ с
областью определения $X_{1}$. Оператор $J$ определяется парой $X$ однозначно; он
называется порождающим для $X$.

Обозначим через $\mathcal{B}$ множество всех измеримых по Борелю функций
$\psi:(0,\infty)\rightarrow(0,\infty)$, которые отделены от нуля на каждом множестве
$[r,\infty)$ и ограниченны на каждом отрезке $[a,b]$, где $r>0$ и $0<a<b<\infty$.

Пусть $\psi\in\mathcal{B}$. В пространстве $X_{0}$ определен, как функция от $J$,
оператор $\psi(J)$. Обозначим через $[X_{0},X_{1}]_\psi$ или, короче, $X_{\psi}$
область определения оператора $\psi(J)$, наделенную скалярным произведением $(w_1,
w_2)_{X_\psi}:=(\psi(J)w_1,\psi(J)w_2)_{X_0}$ и соответствующей нормой
$\|w\|_{X_\psi} = (w,w)_{X_\psi}^{1/2}$. Пространство $X_\psi$ гильбертово и
сепарабельно, причем выполняется непрерывное и плотное вложение $X_\psi
\hookrightarrow X_0$.

Функцию $\psi\in\mathcal{B}$ называем интерполяционным параметром, если для
произвольных допустимых пар $X=[X_0, X_1]$, $Y=[Y_0, Y_1]$ гильбертовых пространств
и для любого линейного отображения $T$, заданного на $X_0$, выполняется следующее.
Если при каждом $j\in\{0,1\}$ сужение отображения $T$ на пространство $X_{j}$
является ограниченным оператором $T:X_{j}\rightarrow Y_{j}$, то и сужение
отображения $T$ на пространство $X_\psi$ является ограниченным оператором
$T:X_{\psi}\rightarrow Y_{\psi}$. Тогда будем говорить, что пространство $X_\psi$
получено интерполяцией с функциональным параметром $\psi$ пары $X$.

Известно, что функция $\psi\in\mathcal{B}$ является интерполяционным параметром
тогда и только тогда, когда она псевдовогнута в окрестности бесконечности, т.~е.
$\psi(t)\asymp\psi_{1}(t)$ при $t\gg1$ для некоторой положительной вогнутой функции
$\psi_{1}(t)$. (Как обычно, $\psi\asymp\psi_{1}$ обозначает ограниченность обоих
отношений $\psi/\psi_{1}$ и $\psi_{1}/\psi$ на указанном множестве.)

В случае, когда допустимая пара состоит из соболевских пространств, нам понадобится
следующий факт \cite{MikhailetsMurach10} (п. 2.4.2, теорема 2.19).

\textbf{Предложение 1.} \it Пусть заданы функция $\varphi\in\mathrm{RO}$ и
вещественные числа $s_0$, $s_1$ такие, что $s_0<\sigma_0(\varphi)$ и
$s_1>\sigma_1(\varphi)$. Положим
\begin{equation}\label{f10a}
\psi(t):=
\begin{cases}
\;t^{{-s_0}/{(s_1-s_0)}}\,
\varphi(t^{1/{(s_1-s_0)}})&\text{при}\quad t\geq1, \\
\;\varphi(1)&\text{при}\quad0<t<1.
\end{cases}
\end{equation}
Тогда функция $\psi\in\mathcal{B}$ является интерполяционным параметром и
\begin{equation}\label{f10b}
[H^{(s_0)},H^{(s_1)}]_{\psi}=H^{\varphi}
\end{equation}
с равенством норм. \rm

Отметим также \cite{MikhailetsMurach10} (п. 2.4.2), что используемый нами класс
гильбертовых пространств $\{H^{\varphi}:\varphi\in\mathrm{RO}\}$ замкнут
относительно интерполяции с функциональным параметром. Более того, он совпадает (с
точностью до эквивалентности норм) с классом всех гильбертовых пространств,
интерполяционных для пар соболевских пространств $[H^{(s_0)},H^{(s_1)}]$, где
$s_0,s_1\in\mathbb{R}$ и $s_0<s_1$. Напомним, что свойство (гильбертового)
пространства $H$ быть интерполяционным для допустимой пары $X=[X_0, X_1]$ означает
следующее: а) выполняются непрерывные вложения $X_1\hookrightarrow H\hookrightarrow
X_0$, б) всякий линейный оператор, ограниченный на каждом из пространств $X_0$ и
$X_1$, является ограниченным и на $X$.

При интерполяции пространств наследуется не только ограниченность, но и нетеровость
линейных операторов при некоторых дополнительных условиях. Сформулируем этот
результат применительно к рассмотренному нами методу интерполяции
\cite{MikhailetsMurach10} (п. 1.1.7, теорема 1.7).

\textbf{Предложение 2.} \it Пусть $X=[X_0,X_1]$ и $Y=[Y_0,Y_1]$~--- допустимые пары
гильбертовых пространств. Пусть, кроме того, на $X_0$ задано линейное отображение
$T$ такое, что его сужения на пространства $X_j$, где $j=0,1$, являются
ограниченными нетеровыми операторами $T:X_j\rightarrow Y_j$, имеющими общее ядро и
одинаковый индекс. Тогда для произвольного интерполяционного параметра
$\psi\in\mathcal{B}$ ограниченный оператор $T:X_\psi\rightarrow Y_\psi$ нетеров с
теми же ядром и индексом, а его область значений равна $Y_\psi\cap T(X_0)$. \rm

В доказательствах нам придется интерполировать ортогональные суммы гильбертовых
пространств. Для этого будет полезен следующий факт \cite{MikhailetsMurach10} (п.
1.1.5, теорема 1.5).

\textbf{Предложение 3.} \it Пусть задано конечное число допустимых пар
$[X_{0}^{(k)},X_{1}^{(k)}]$ гильбертовых пространств, где $k=1,\ldots,p$. Тогда для
любого $\psi\in\mathcal{B}$ справедливо
$$
\biggl[\,\bigoplus_{k=1}^{p}X_{0}^{(k)},\,\bigoplus_{k=1}^{p}X_{1}^{(k)}\biggr]_{\psi}=\,
\bigoplus_{k=1}^{p}\bigl[X_{0}^{(k)},\,X_{1}^{(k)}\bigr]_{\psi}
$$
с равенством норм. \rm

При доказательстве теорем 1 и 2 мы воспользуемся тем важным фактом, что равномерно
эллиптический дифференциальный оператор $A$ имеет параметрикс, т.~е. матричный
псевдодифференциальный оператор (ПДО), обратный к $A$ с точностью до ПДО порядка
$-\infty$. Напомним необходимые нам факты, относящиеся к ПДО и параметриксам (см.,
например, \cite{Agranovich94} (п. 1.1, 1.9, 3.2)).

Обозначим через $\Psi^r$, где $r\in\mathbb{R}$, множество всех ПДО $G$ в
$\mathbb{R}^n$ (не обязательно классических) таких, что их символ $g(x,\xi)$
бесконечно дифференцируемый в $\mathbb{R}^{2n}$ и удовлетворяет следующему условию:
для любых мультииндексов $\alpha$ и $\beta$ существует число $c_{\alpha,\beta}> 0$,
при котором
$$
|\,D_x^{\alpha}D_\xi^{\beta}\emph{g}(x,\xi)|\leq
c_{\alpha,\beta}\langle\xi\rangle^{r-|\beta|}\quad\mbox{для любых}\quad x,
\xi\in\mathbb{R}^{n}.
$$
Число $r$ называется (формальным) порядком ПДО $G$. Положим
$\Psi^{-\infty}:=\bigcap_{r\in\mathbb{R}}\Psi^r$.

\textbf{Предложение 4.} \it Существует матричный ПДО $B=(B_{k,j})_{k,j=1}^p$ такой,
что все $B_{k,j}\in\Psi^{-m_k-l_j}$ и
\begin{equation}\label{f11}
BA=I+T_1,\quad  AB=I+T_2,
\end{equation}
где   $T_1= (T_1^{j,k})_{j,k=1}^p$ и $T_2=(T_2^{k,j})_{k,j=1}^p$~--- некоторые
матричные ПДО, состоящие из элементов класса $\Psi^{-\infty}$,   а I~---
тождественный оператор в $S'$. \rm

Всякий ПДО класса $\Psi^r$ является непрерывным оператором в пространстве
$\mathcal{S}'$. Следующая лемма уточняет этот факт применительно к пространствам
Хермандера.

\textbf{Лемма 1.} \it Для ПДО $G\in\Psi^r$ сужение линейного отображения
$u\rightarrow Gu$, $u\in\mathcal{S}'$, на пространство $H^\varphi$ является
ограниченным оператором
\begin{equation}\label{f12}
G:\,H^{\varphi}\rightarrow H^{\varphi\rho^{-r}}\quad\mbox{для
любых}\quad\varphi\in\mathrm{RO}.
\end{equation}\rm

\textbf{\textit{Доказательство.}} В случае соболевских пространств этот факт
известен \cite{Agranovich94} (п. 1.1, теорема 1.1.2). Отсюда выведем ограниченность
оператора \eqref{f12} с помощью интерполяции с функциональным параметром.

Пусть $\varphi\in\mathrm{RO}$. Выберем числа $s_{0}<\sigma_{0}(\varphi)$ и
$s_{1}>\sigma_{1}(\varphi)$. Рассмотрим линейные ограниченные операторы
\begin{equation}\label{f12a}
G:\,H^{(s_{j})}\rightarrow H^{(s_{j}-r)}\quad\mbox{для}\quad j=0,1,
\end{equation}
действующие в пространствах Соболева. Определим $\psi$ по формуле $\eqref{f10a}$;
согласно предложению~1, функция $\psi$--- интерполяционный параметр. Потому из
ограниченности операторов \eqref{f12a} следует, что сужение отображения $G$ на
пространство $[H^{(s_0)},H^{(s_1)}]_\psi$ является ограниченным оператором
\begin{equation}\label{f12b}
G:\bigl[H^{(s_0)},H^{(s_1)}\bigr]_\psi\rightarrow
\bigl[H^{(s_0-r)},H^{(s_1-r)}\bigr]_\psi.
\end{equation}

В силу предложения~1 выполняются равенства $\eqref{f10b}$ и
$$
\bigl[H^{(s_0-r)},H^{(s_1-r)}\bigr]_\psi=H^{\varphi\rho^{-r}}.
$$
Заметим, что второе из них верно, поскольку $s_{0}-r<\sigma_{0}(\varphi\rho^{-r})$,
$s_{1}-r>\sigma_{1}(\varphi\rho^{-r})$, а функциональный параметр $\psi$
удовлетворяет соотношению $\eqref{f10a}$, если в нем заменить $s_0$ на $s_0-r$,
$s_1$ на $s_1-r$ и $\varphi$ на $\varphi\rho^{-r}$. Следовательно, ограниченность
оператора $\eqref{f12b}$ означает ограниченность оператора $\eqref{f12}$.

Лемма~1 доказана.

В силу леммы~1 оператор \eqref{f2} ограничен, поскольку каждый дифференциальный
оператор $A_{j,k}(x,D)$ принадлежит классу $\Psi^{l_j+m_k}$.

Для доказательства теоремы~3 нам понадобится следующий изотропный вариант теоремы
вложения Хермандера.

\textbf{Лемма 2.} \it Пусть заданы целое число $\lambda\geq0$ и функция
$\omega\in\mathrm{RO}$. Тогда условие
\begin{equation}\label{f12c}
\int\limits_{1}^{\infty}\,t^{2\lambda+n-1}\,\omega^{-2}(t)\,dt<\infty
\end{equation}
равносильно вложению $H^{\omega}\subset C^{\lambda}_{\mathrm{b}}$, и это вложение
непрерывно. \rm

\textbf{\textit{Доказательство.}} Теорема вложения Хермандера \cite{Hermander63} (п.
2.2, теорема 2.2.7) утверждает в гильбертовом случае, что
$$
\int\limits_{\mathbb{R}^{n}}\,\langle\xi\rangle^{2\lambda}\,k^{-2}(\xi)\,d\xi<\infty\;
\Leftrightarrow\;B_{2,k}\subset C^{\lambda}_{\mathrm{b}}.
$$
Здесь, напомним, $B_{2,k}$~--- пространство Хермандера, параметризуемое весовой
функцией $k(\xi)$ от $n$ переменных. Если эта функция радиальна:
$k(\xi)=\omega(\langle\xi\rangle)$, то соответствующее ей пространство
$B_{2,k}=H^{\omega}$ изотропно и
\begin{equation}\label{f12d}
\int\limits_{\mathbb{R}^{n}}\,\langle\xi\rangle^{2\lambda}\,
\omega^{-2}(\langle\xi\rangle)\,d\xi<\infty\;\Leftrightarrow\;H^{\omega}\subset
C^{\lambda}_{\mathrm{b}}.
\end{equation}
Покажем, что левое условие в \eqref{f12d} эквивалентно \eqref{f12c}.

Переходя к сферическим координатам, где $r:=|\xi|$, и затем делая замену
$t=\sqrt{1+r^{2}}$, получаем:
\begin{gather*}
\int\limits_{\mathbb{R}^{n}}\,\langle\xi\rangle^{2\lambda}\,
\omega^{-2}(\langle\xi\rangle)\,
d\xi=c\,\int\limits_{0}^{\infty}\,(1+r^{2})^{\lambda}\,\omega^{-2}(\sqrt{1+r^{2}})\,
r^{n-1}\,dr=\\=
c\,\int\limits_{1}^{\infty}\,t^{2\lambda+1}\,(t^{2}-1)^{n/2-1}\,\omega^{-2}(t)\,dt=
A+c\,\int\limits_{2}^{\infty}\,t^{2\lambda+1}\,(t^{2}-1)^{n/2-1}\,\omega^{-2}(t)\,dt.
\end{gather*}
Здесь $c:=nV_{1}$, где $V_{1}$~--- объем единичного шара в $\mathbb{R}^{n}$, а
$$
A:=c\,\int\limits_{1}^{2}\,t^{2\lambda+1}\,(t^{2}-1)^{n/2-1}\,\omega^{-2}(t)\,dt<\infty,
$$
поскольку $\omega\asymp1$ на $[1,2]$ и $n/2-1>-1$. Следовательно,
\begin{gather*}
\int\limits_{\mathbb{R}^{n}}\,\langle\xi\rangle^{2\lambda}\,
\omega^{-2}(\langle\xi\rangle)\,d\xi<\infty\;\Leftrightarrow\;
\int\limits_{2}^{\infty}\,t^{2\lambda+1}\,(t^{2}-1)^{n/2-1}\,\omega^{-2}(t)\,dt<\infty
\;\Leftrightarrow\\ \Leftrightarrow\;
\int\limits_{2}^{\infty}\,t^{2\lambda+n-1}\,\omega^{-2}(t)\,dt<\infty\;
\Leftrightarrow\;\eqref{f12c}.
\end{gather*}
Отсюда в силу \eqref{f12d} делаем вывод, что \eqref{f12c} эквивалентно вложению
$H^{\omega}\subset C^{\lambda}_{\mathrm{b}}$. Оно непрерывно, поскольку банаховы
пространства $H^{\omega}$ и $C^{\lambda}_{\mathrm{b}}$ непрерывно вложены в
некоторое хаусдорфово пространство, например в $\mathcal{S}'$.

Лемма~2 доказана.

В связи с ней отметим следующее. Если $H^{\omega}=H^{(s)}$~--- пространство Соболева
порядка $s$, т.~е. $\omega(t)=t^{s}$ при $t\geq1$, то условие \eqref{f12c}
равносильно неравенству $s>\lambda+n/2$, и мы приходим к теореме вложения Соболева.

\textbf{5. Доказательство основных результатов.} Докажем теоремы 1~-- 4.

\textbf{\textit{Доказательство теоремы}~1\textit{.}} Обозначим через
$\|\cdot\|_{\varphi}'$, $\|\cdot\|_{\varphi}''$ и $\|\cdot\|_ {\varphi,\sigma}'$
соответственно нормы в пространствах
$$
\bigoplus_{k=1}^{p}H^{\varphi\rho^{m_k}},\quad
\bigoplus_{j=1}^{p}H^{\varphi\rho^{-l_j}}\quad\mbox{и}\quad
\bigoplus_{k=1}^{p}H^{\varphi\rho^{m_k-\sigma}}.
$$

Пусть вектор-функции \eqref{f3} удовлетворяют уравнению $Au=f$ в $\mathbb{R}^{n}$. В
силу первого равенства в $\eqref{f11}$ имеем: $u=Bf-T_1u$. Отсюда вытекает оценка
\eqref{f4}:
$$
\|u\|_{\varphi}'=\|Bf-T_{1}u\|_{\varphi}'\leq\|Bf\|_{\varphi}'+\|T_{1}u\|_{\varphi}'
\leq c\,\|f\|_{\varphi}'+c\,\|u\|_{\varphi, \sigma}'.
$$
Здесь с~--- максимум норм операторов
\begin{gather}\label{f13}
B:\,\bigoplus_{j=1}^{p}H^{\varphi\rho^{-l_j}}\rightarrow
\bigoplus_{k=1}^{p}H^{\varphi\rho^{m_k}},\\
T_1\,:\bigoplus_{k=1}^{p}H^{\varphi\rho^{m_k-\sigma}}\rightarrow
\bigoplus_{k=1}^{p}H^{\varphi\rho^{m_k}},\notag
\end{gather}
ограниченных в силу предложения~4 и леммы~1.

Теорема 1 доказана.

\textbf{\textit{Доказательство теоремы}~2\textit{.}} Сначала рассмотрим случай,
когда $V=\mathbb{R}^{n}$. По условию, $Au=f$ в $\mathbb{R}^{n}$, где
$f\in\bigoplus_{j=1}^{p}H^{\varphi\rho^{-l_j}}$. Воспользовавшись первым равенством
в $\eqref{f11}$, запишем $u=Bf-T_1u$. Здесь
$Bf\in\bigoplus_{k=1}^{p}H^{\varphi\rho^{m_k}}$ в силу \eqref{f13} и
$T_1u\in(H^{\infty})^{p}$ ввиду предложения~4. Следовательно,
$u\in\bigoplus_{k=1}^{p}H^{\varphi\rho^{m_k}}$, что и требовалось доказать в случае
$V=\mathbb{R}^{n}$.

Рассмотрим теперь случай, когда $V\neq\mathbb{R}^{n}$. Произвольно выберем функцию
$\chi\in C^{\infty}_{\mathrm{b}}$ такую, что
\begin{equation}\label{14a}
\mathrm{supp}\,\chi\subset V\quad\mbox{и}\quad
\mathrm{dist}(\mathrm{supp}\,\chi,\partial V)>0.
\end{equation}
Для нее существует функция $\eta\in C^{\infty}$ такая, что
\begin{equation}\label{14b}
\mathrm{supp}\,\eta\subset V,\;\;\mathrm{dist}(\mathrm{supp}\,\eta,\,\partial
V)>0,\;\;\eta=1\;\mbox{в окрестности}\;\mathrm{supp}\,\chi.
\end{equation}
Действительно, можно определить указанную функцию с помощью операции свертки по
формуле $\eta:=\chi_{2\varepsilon}\ast\omega_{\varepsilon}$, где
$\varepsilon:=\mathrm{dist}(\mathrm{supp}\,\chi,\partial V)/4$,
$\chi_{2\varepsilon}$ --- индикатор $2\varepsilon$-окрестности множества
$\mathrm{supp}\,\chi$, а функция $\omega_{\varepsilon}\in C^{\infty}_{0}$
удовлетворяет условиям
$$
\omega_{\varepsilon}\geq0,\;\;
\mathrm{supp}\,\omega_{\varepsilon}\subset\{x\in\mathbb{R}^{n}:\|x\|\leq
\varepsilon\},\;\;\int\limits_{\mathbb{R}^{n}}\omega_{\varepsilon}(x)\,dx=1.
$$
Непосредственно проверяется, что такая функция $\eta$ принадлежит классу
$C^{\infty}_{\mathrm{b}}$ и имеет следующее свойство: $\eta\equiv1$ в
$\varepsilon$-окрестности множества $\mathrm{supp}\,\chi$ и $\eta\equiv0$ вне
$3\varepsilon$-окрестности этого же множества, т.~е. $\eta$ удовлетворяет условиям
\eqref{14b}.

На основании первого равенства в \eqref{f11} можем записать
\begin{equation}
\begin{gathered}\label{14c}
\chi u=\chi BAu-\chi T_{1}u=\chi B\eta Au+\chi B(1-\eta)Au-\chi T_{1}u.
\end{gathered}
\end{equation}
Так как $Au=f$ на множестве $V$, то $\eta Au=\eta f$ в $\mathbb{R}^{n}$, где $\eta
f\in\bigoplus_{j=1}^{p}H^{\varphi\rho^{-l_j}}$ по условию теоремы. Следовательно, в
силу $\eqref{f13}$ имеем:
$$
\chi B\eta Au=\chi B\eta f\in\bigoplus_{k=1}^p H^{\varphi\rho^{m_k}}.
$$
Кроме того, поскольку матричные ПДО $\chi B(1-\eta)$, где $1-\eta=0$ в окрестности
$\mathrm{supp}\,\chi$, и $T_1$ состоят из элементов класса $\Psi^{-\infty}$, то
вектор-функции $\chi B(1-\eta)Au$ и $T_{1}u$ принадлежат пространству
$(H^\infty)^{p}$. Поэтому в силу \eqref{14c} получаем, что $\chi u
\in\bigoplus_{k=1}^{p}H^{\varphi\rho^{m_k}}$ для любой функции $\chi\in
C^{\infty}_{\mathrm{b}}$, удовлетворяющей условию \eqref{14a}. Иными словами,
$u_{k}\in H_{\mathrm{int}}^{\varphi\rho^{m_k}}(V)$ для всех $k=1,\ldots,p$.

Теорема~2 доказана.

\textbf{\textit{Доказательство теоремы}~3\textit{.}} Сначала рассмотрим случай,
когда $V=\mathbb{R}^{n}$. В силу теоремы~2 имеем: $u_k\in H^{\varphi\rho^{m_k}}$.
Отсюда на основании леммы~2, где $\omega:=\varphi\rho^{m_k}$, и условия $\eqref
{f6}$ получаем включение $u_k\in C_{\mathrm{b}}^\lambda$, что и требовалось доказать
в этом случае.

Предположим теперь, что $V\neq\mathbb{R}^{n}$. В силу теоремы~2 имеем: $u_k\in
H_{\mathrm{int}}^{\varphi\rho^{m_k}}(V)$. Пусть функция $\eta\in
C^{\infty}_{\mathrm{b}}$ удовлетворяет следующим условиям:
$\mathrm{supp}\,\eta\subset V$, $\mathrm{dist}(\mathrm{supp}\,\eta,\,\partial V)>0$
и $\eta=1$ в окрестности множества $V_{0}\subset V$ таком, что
$\mathrm{dist}(V_{0},\partial V)>0$. Эта функция строится так же как и в
доказательстве теоремы~2, если заменить в нем множество $\mathrm{supp}\,\chi$ на
$V_{0}$. На основании леммы~2, где $\omega:=\varphi\rho^{m_k}$, и условия $\eqref
{f6}$ имеем: $\eta u_k\in H^{\varphi\rho^{m_k}}\subset C_{\mathrm{b}}^\lambda$.
Отсюда следует, что все частные производные функции $u_k$ до порядка $\lambda$
включительно непрерывны и ограничены в некоторой окрестности множества $V_0$. Тогда
эти производные непрерывны и на множестве $V$, поскольку можно взять  $V_0:=\{x_0\}$
для любой точки $x_0\in V$.

Теорема~3 доказана.

\textbf{\textit{Доказательство теоремы}~4\textit{.}} В соболевском случае
$\varphi=\varrho^{s}$, где произвольно выбрано $s\in\mathbb{R}$, эта теорема
известна (см., например, \cite{Rabier12} (теорема 4.2)). Отсюда выведем ее для
любого $\varphi\in\mathrm{RO}$ с помощью интерполяции с функциональным параметром.

Выберем числа $s_{0}<\sigma_{0}(\varphi)$ и $s_{1}>\sigma_{1}(\varphi)$. Рассмотрим
ограниченные нетеровы операторы
\begin{equation}\label{f16}
A:\,\bigoplus_{k=1}^{p}H^{(s_r+m_k)}\rightarrow
\bigoplus_{j=1}^{p}H^{(s_r-l_j)}\quad\mbox{для}\quad r=0,1,
\end{equation}
действующие в пространствах Соболева. Эти операторы имеют общее ядро $N$, одинаковый
индекс, равный $\dim N-\dim N^{+}$, и области значений
\begin{equation}\label{f17}
A\biggl(\,\bigoplus_{k=1}^{p}H^{(s_r+m_k)}\biggr)=
\biggl\{f\in\bigoplus_{j=1}^{p}H^{(s_r-l_j)}:\,(f,v)_{\mathbb{R}^{n}}=0\;\;\mbox{для
всех}\;\;v\in N^{+}\biggr\}.
\end{equation}

Определим интерполяционный параметр $\psi$ по формуле $\eqref{f10a}$. Согласно
предложению~2 нетеровость операторов \eqref{f16} влечет за собой нетеровость
ограниченного оператора
\begin{equation}\label{f18}
A:\,
\biggl[\,\bigoplus_{k=1}^{p}H^{(s_0+m_k)},\bigoplus_{k=1}^{p}H^{(s_1+m_k)}\biggr]_\psi
\rightarrow
\biggl[\,\bigoplus_{j=1}^{p}H^{(s_0-l_j)},\bigoplus_{j=1}^{p}H^{(s_1-l_j)}\biggr]_\psi.
\end{equation}
Здесь в силу предложений 3 и 1 имеем:
\begin{gather}\label{f19}
\biggl[\,\bigoplus_{k=1}^{p}H^{(s_0+m_k)},\bigoplus_{k=1}^{p}H^{(s_1+m_k)}\biggr]_\psi=
\,\bigoplus_{k=1}^{p}\bigl[H^{(s_0+m_k)},H^{(s_1+m_k)}\bigr]_\psi=
\,\bigoplus_{k=1}^{p}H^{\varphi\rho^{m_k}},\\
\biggl[\,\bigoplus_{j=1}^{p}H^{(s_0-l_j)},\bigoplus_{j=1}^{p}H^{(s_1-l_j)}\biggr]_\psi=
\,\bigoplus_{j=1}^{p}\bigl[H^{(s_0-l_j)},H^{(s_1-l_j)}\bigr]_\psi=
\,\bigoplus_{j=1}^{p}H^{\varphi\rho^{-l_j}}. \label{f20}
\end{gather}
Уточним, что \eqref{f19} верно на основании предложения~1, поскольку
$s_{0}+m_k<\sigma_{0}(\varphi\rho^{m_k})$,
$s_{1}+m_k>\sigma_{1}(\varphi\rho^{m_k})$, а функциональный параметр $\psi$
удовлетворяет соотношению $\eqref{f10a}$, если в нем заменить $s_0$ на $s_0+m_k$,
$s_1$ на $s_1+m_k$ и $\varphi$ на $\varphi\rho^{m_k}$. Аналогично, \eqref{f20}
верно, поскольку $s_{0}-l_j<\sigma_{0}(\varphi\rho^{-l_j})$,
$s_{1}-l_j>\sigma_{1}(\varphi\rho^{-l_j})$, а $\psi$ удовлетворяет соотношению
$\eqref{f10a}$, если в нем заменить $s_0$ на $s_0-l_j$, $s_1$ на $s_1-l_j$ и
$\varphi$ на $\varphi\rho^{-l_j}$.

Таким образом, \eqref{f2}~--- это нетеров оператор \eqref{f18}. В силу
предложения~2, индекс оператора \eqref{f2} равен $\dim N-\dim N^{+}$, а область
значений равна
$$
\bigoplus_{j=1}^{p}H^{\varphi\rho^{-l_j}}\cap
A\biggl(\,\bigoplus_{k=1}^{p}H^{(s_0+m_k)}\biggr)
$$
и совпадает с \eqref{f7} ввиду \eqref{f17}.

Теорема 4 доказана.

\vspace{1cm}

Институт математики Национальной академии наук Украины,

ул. Терещенкивська, 3, 01601 Киев-4, Украина

\vspace{0.5cm}

Institute of Mathematics of National Academy of Sciences of Ukraine,

3, Tereshchenkivska Str., 01601 Kyiv-4, Ukraine

\vspace{0.5cm}

\verb"zinchenko@imath.kiev.ua"

\verb"murach@imath.kiev.ua"

\end{document}